\documentclass[a4paper,11pt]{article}
\usepackage{graphicx}
\usepackage{epstopdf}
\usepackage{amsthm}
\usepackage[utf8]{inputenc}
\usepackage[english]{babel}
\newtheorem{theorem}{Theorem}[section]

\theoremstyle{definition}

\usepackage{amssymb, amsmath}
\title{Certain subclass of harmonic univalent functions defined by $q$-differential operator}
\date{}
\begin{document}
\numberwithin{equation}{section}
\maketitle
\date{}
\begin{center}
\author{{\bf{G. M. Birajdar}}\vspace{.31cm}\\
School of Mathematics \& Statistics,\\
Dr. Vishwanath Karad MIT World Peace University,\\
Pune (M.S) India 411038\\
	Email: gajanan.birajdar@mitwpu.edu.in}\vspace{.5cm}\\
 \author{{\bf{N. D. Sangle}}\vspace{0.31cm}\\
 	Department of Mathematics,\\
 	D. Y. Patil College of Engineering \& Technology, \\
 	Kasaba Bawada,
 	Kolhapur, (M.S.), India 416006\\
 	Email: navneet\_sangle@rediffmail.com}
 \end{center}
\vspace{1cm}
\abstract{}
In this paper, we define certain subclass of harmonic univalent function in the unit disc $U = \left\{ {z \in C:\left| z \right| < 1} \right\}$
by using $q$- differential operator. Also we obtain coefficient inequalities,  growth and distortion theorems for this subclass. \\
 

\section{Introduction} 
Clunie and Sheil-Small \cite{B1} investigated the class $S_H$ as well as its geometric subclasses and established some coefficient
bounds. Since then, there have been several related papers on $S_H$ and its subclasses. In fact, by introducing new subclasses, Silverman \cite{B11}, Silverman and Silvia \cite{B12}, Jahangiri \cite{B7}, Sangle and Yadav \cite{B8}, Dixit and Porwal \cite{B4}, Singh and Porwal \cite{B13} and Ravindar et.al \cite{B14} etc. presented a systematic and unified study of harmonic univalent functions.\\
The concepts of q-calculus has many applications in subfields of science, some of them are q-difference equations and geometric function theory.
Motivated by the research work done by Jahangiri \cite{B2,B3}, Joshi and Sangle \cite{B5,B9}, Purohit et al.\cite{B6}, we define some subclasses of harmonic mappings using the Salagean $q$-differential operator.\\
Also, we determine extreme points and coefficient estimates of $ S^{q}_{H}(m,\alpha,u)$ and $\overline{S}^{q}_{H}(m,\alpha,u)$.\\
Let A be family of analytic functions in unit disk $U$ and 
${A^0}$ be the class of all normalized analytic functions. For $0 < q < 1$ and for positive integer $u$, the $q$-integer number is denoted by ${\left[u \right]_q}$ and also it is written as 
\begin{equation}
{\left[ u \right]_q} = \frac{{1 - {q^u}}}{{1 - q}}=\sum\limits_{k \ge 0}^{u - 1} {{q^k}}.
\end{equation}
By making use of differential calculus, we can check that 
$\mathop {\lim }\limits_{q \to {1^ -}} {\left[u\right]_q} = u$ \\
For $h \in A$, the $q$-difference operator \cite{B6} is specified as 
\begin{align}
{\partial _q}h(z) = \frac{{h(z) - h(qz)}}{{\left( {1 - q} \right)z}}
\end{align}
where $\mathop {\lim }\limits_{q \to {1^ -}} {\partial _q}h(z) = h'(z)$.\\
Let the functions $h \in A$ be of the form \\
\begin{align}
h(z) = z + \sum\limits_{u \ge 2}^\infty  {{a_u}\,{z^u}}.
\end{align}
J. M. Jahangiri \cite{B3} defined the Salagean $q$-differential operator for the above functions $h$ as 
\begin{equation*}
D_q^0h(z) = h(z)
\end{equation*}
\begin{equation*}
D_q^1h(z) = z{\partial _q}h(z) = \frac{{h(z) - h(qz)}}{{\left( {1 - q} \right)z}},...
\end{equation*}
\begin{equation}
D_q^mh(z) = z{\partial _q}D_q^{m - 1}h(z) = h(z)\, * \left( {z + \sum\limits_{u \ge 2}^\infty  {\left[ u \right]_q^m\,{z^u}} } \right)  = z + \sum\limits_{u \ge 2}^\infty  {\left[ u \right]_q^m\,{a_u}\,{z^u}}
\end{equation}
where $m$ is a positive integer. The operator  $D^m_q$   is called Salagean q-differential operator.\\
The complex-valued harmonic functions can be written as $f = h + \overline g$ in  where $h$ and $g$ have the following power series expansions
\begin{equation}
h(z) = z + \sum\limits_{u \ge 2}^\infty  {{a_u}\,{z^u}},\quad
g(z) = \sum\limits_{u \ge 1}^\infty  {{b_u}\,{z^u}} ,
\left| {{b_1}} \right| < 1.
\end{equation}
Clunie and Sheil-Small \cite{B1} defined the function of form $f = h + \overline g$ that are locally univalent, sense-preserving and harmonic in $U$. A sufficient condition for the harmonic functions $f$ to be univalent in $U$ is that
$\left| {h'(z)} \right| \ge \left| {g'(z)} \right|$ in $U$.\\ 
J. M. Jahangiri \cite{B3} defined the Salagean $q$-differential operator for the harmonic functions $f$ by
\begin{align}
D_q^mf(z) = D_q^mh(z) + {\left({- 1} \right)^m}\,\overline {D_q^m g(z)}
\end{align}\\
where $D_q^m$ is defined by (1.4).\\
 Now, for $0\le \alpha<1 ,m \in N_{0} $ and $z \in U$, suppose that $ S^{q}_{H}(m,\alpha,u)$ denote the family of harmonic univalent function $f$ of the form $f=h+\overline g$ such that 
 \begin{align}
 Re\bigg\{\frac{D^{m}_q h(z)+D^{m}_q g(z)}{z}\bigg\}>\alpha
 \end{align}
 where $D_q^mf(z)$ is defined by J. M. Jahangiri \cite{B3}.\\
Further let the subclass $\overline{S}^{q}_{H}(m,\alpha,u)$ consisting harmonic functions  $f=h+\overline{g} $ \ in \  $\overline{S}^{q}_{H}(m,\alpha,u)$ so that $h$ and $g$ are of the form 
\begin{align}
h(z)=z-\sum_{u\ge2}^{\infty}{\left|a_{u}\right|}{z^{u}} \quad  and \quad   g(z)=\sum_{u\ge1}^{\infty}{\left|b_{u}\right|}{z^{u}}.
\end{align}
\section{Main Results} 

\begin{theorem}
Let the function $f=h+\overline{g}$ be such that $h$ and $g$ are given by (1.5), \\
Furthermore 	
\begin{equation}
\sum\limits_{u \ge 2}^\infty  {\left[ u \right]} _q^m\,\left| {{a_u}} \right|\, + \sum\limits_{u \ge 1}^\infty  {\left[ u \right]} _q^m\,\left| {{b_u}} \right|\, \le \left( {1 - \alpha } \right)
\end{equation}
where $0\le\alpha<1$ and $ m \in N_{0}$. Then $f$ is harmonic univalent, sense-preserving in $U$ and $f  \in  S^{q}_{H}(m,\alpha,u)$.
\end{theorem}
\textbf{Proof:}
If \ $ z_1 \neq z_2 $  then,
\begin{align*}
\left| {\frac{{f({z_1}) - f({z_2})}}{{h({z_1}) - h({z_2})}}} \right| &\ge1-\left| {\frac{{g({z_1}) - g({z_2})}}{{h({z_1}) - h({z_2})}}} \right|\\
&=1-\left|\frac{{\sum\limits_{u \ge 1}^{ \infty } {{b_u}} ({z_1}^u - {z_2}^u)}}{z_1 -z_2 +\sum\limits_{k \ge 2}^{\infty } {{a_u}} ({z_1}^u - {z_2}^u)}\right|\\
&>1-\frac{\sum\limits_{u \ge 1}^{ \infty} {u \left|b_u\right|}}{1-\sum\limits_{u \ge 2}^{\infty } {u \left|{a_u}\right|}}\\
&\ge 1 - \frac{{\sum\limits_{u \ge 1}^\infty  {\frac{{\left[ u \right]_q^m}}{{1 - \alpha }}\,\left| {{b_u}} \right|} }}{{1 - \sum\limits_{u \ge 2}^\infty  {\frac{{\left[ u \right]_q^m}}{{1 - \alpha }}\,\left| {{a_u}} \right|} }}\\
&\ge0.
\end{align*}
Hence $f$ is univalent in $U$.\\
$f$ is sense-preserving in $U$. This is because
\begin{align*}
\left| {h'(z)}\right|&\ge1-\sum\limits_{u\ge 2}^\infty{u\,\left|{{a_u}} \right|} {\left| z \right|^{u - 1}}\\
&>1-\sum\limits_{u\ge 2}^\infty{u\,\left| {{a_u}} \right|}\\
&\ge{1 - \sum\limits_{u \ge 2}^\infty{\frac{{\left[ u \right]_q^m}}{{1-\alpha }}\,\left| {{a_u}} \right|} }\\
&\ge{\sum\limits_{u\ge1}^\infty{\frac{{\left[u\right]_q^m}}{{1-\alpha }}\,\left| {{b_u}} \right|} }\\
&\ge \sum\limits_{u \ge 1}^\infty  {u\,\,\left| {{b_u}} \right|} {\left| z \right|^{u - 1}}\\
&\ge \left| {g'(z)} \right|.
\end{align*}
Now, we show that  $f \in S^{q}_{H}(m,\alpha,u)$.Using the fact that 
${\mathop{\rm Re}\nolimits} \left( w \right) > \alpha$ if and only if 
$\left| {1 - \alpha  + w} \right| > \left| {1 + \alpha  - w} \right|$.
it suffices to show that
\begin{align}
\left| {(1 - \alpha)  + \frac{{D_q^mh(z) + D_q^mg(z)}}{z}} \right| - \left| {(1 + \alpha)  - \frac{{D_q^mh(z) + D_q^mg(z)}}{z}} \right| >0
\end{align}
Substituting for $D_q^mh(z)$ and $D_q^mg(z)$ in (2.2), we obtain
\begin{align*}
&= \left| {(2 - \alpha ) + \sum\limits_{u \ge 2}^\infty  {\left[ u \right]_q^m\,{a_u}\,{z^{u - 1}}}  + \sum\limits_{u \ge 1}^\infty  {\left[ u \right]_q^m\,{b_u}\,{z^{u - 1}}} } \right|-\left| {  \alpha  -\sum\limits_{u \ge 2}^\infty  {\left[ u \right]_q^m\,{a_u}\,{z^{u - 1}}}  -\sum\limits_{u \ge 1}^\infty  {\left[ u \right]_q^m\,{b_u}\,{z^{u - 1}}} } \right|\\
&\ge 2(1 - \alpha )\left\{ {1 - \sum\limits_{u \ge 2}^\infty  {\frac{{\left[ u \right]_q^m}}{{1 - \alpha }}\,\left| {{a_u}} \right|{{\left| z \right|}^{u - 1}} - \sum\limits_{u \ge 1}^\infty  {\frac{{\left[ u \right]_q^m}}{{1 - \alpha }}\,\left| {{b_u}} \right|{{\left| z \right|}^{u - 1}}} } } \right\}\\
&> 2(1 - \alpha )\left\{ {1 - \sum\limits_{u \ge 2}^\infty  {\frac{{\left[ u \right]_q^m}}{{1 - \alpha }}\,\left| {{a_u}} \right| - \sum\limits_{u \ge 1}^\infty  {\frac{{\left[ u \right]_q^m}}{{1 - \alpha }}\,\left| {{b_u}} \right|} } } \right\}
\end{align*}
The harmonic mappings
\begin{equation*}
f(z) = z + \sum\limits_{u \ge 2}^\infty  {\frac{{1 - \alpha }}{{\left[ u \right]_q^m}}} \,{x_u}\,{z^u} + \sum\limits_{u \ge 1}^\infty  {\frac{{1 - \alpha }}{{\left[ u \right]_q^m}}} \,\overline {{y_u}\,{z^u}},
\end{equation*}\\
where $\sum\limits_{u \ge 2}^\infty  {\left| {{x_u}} \right|} \, + \sum\limits_{u \ge 1}^\infty  {\left| {{y_u}} \right| = 1}$,
show that coefficient bound given by (2.1) is sharp.\\
In the following theorem, it is proved that the condition (2.1) is also necessary for functions $f=h+\overline{g}$ where $h$ and $g$ are
of the form (1.8).
\begin{theorem}
Let $ f = h + \overline g$ be given by (1.8).Then $f \in \overline S _H^q(m,\alpha,u)$ if and only if
\begin{align}
\sum\limits_{u \ge 2}^\infty {\frac{{\left[u\right]_q^m}}{{1-\alpha}}} \,\left| {{a_u}} \right| + \sum\limits_{u \ge 1}^\infty  {\frac{{\left[u \right]_q^m}}{{1 - \alpha }}} \,\left| {{b_u}} \right| \le 1
\end{align}
where $0 \le \alpha  < 1$ and $m \in {N_0}$.
\end{theorem}
\textbf{Proof:}
The if part follows from Theorem 2.1. For the only if part, we show that
$f \in \overline S _H^q(m,\alpha,u)$ if the condition (2.3) holds. We
notice that the condition
\begin{align*}
{\mathop{\rm Re}\nolimits} \left\{ {\frac{{D_q^m\,h(z) + D_q^m\,g(z)}}{z}} \right\} > \alpha
\end{align*}
is equivalent to
\begin{align*}
{\mathop{\rm Re}\nolimits} \left\{ {1 - \sum\limits_{u \ge 2}^\infty  {\left[ u \right]_q^m\,\left| {{a_u}} \right|{{\left| z \right|}^{u - 1}} - \sum\limits_{u \ge 1}^\infty  {\left[ u \right]_q^m\,\left| {{b_u}} \right|{{\left| z \right|}^{u - 1}}} } } \right\} > \alpha. 
\end{align*}
The above required condition must hold for all values of $z$ in $U$. Taking the values of $z$ on the positive real axis, where $0 \le \left| z \right| = r < 1$, we must have
\begin{align*} 1 - \sum\limits_{u \ge 2}^\infty  {\left[ u \right]} _q^m\,\left| {{a_u}} \right|\, - \sum\limits_{u \ge 1}^\infty  {\left[ u \right]} _q^m\,\left| {{b_u}} \right|\, \ge \alpha
\end
{align*}
which is precisely the assertion (2.3).\\
Next, we determine the extreme points of closed convex hulls of class 
$\overline S _H^q(m,\alpha,u)$.
\begin{theorem}
Let $f$ be given by (1.8). Then  $\overline S _H^q(m,\alpha,u)$ if and only if 
\begin{equation*}
f(z) = \sum\limits_{u \ge 1}^\infty  {\left( {{x_u}\,{h_u}(z) + {y_u}\,{g_u}(z)} \right)},
\end{equation*}
where ${h_1}(z) = z$,
\begin{equation*}
{h_k}(z) = z - \frac{{1 - \alpha }}{{\left[ u \right]_q^m}}\,{z^u}\,,\left( {u = \,2,3,4,...} \right)\,,{g_k}(z) = z - \frac{{1 - \alpha }}{{\left[ u \right]_q^m}}\,{\overline z ^u}\,,\left( {u = \,1,2,3,4,...} \right),
\end{equation*}
${x_u} \ge 0\,,{y_u} \ge 0\,,\sum\limits_{u = 1}^\infty  {{x_u} + } {y_u} = 1$.
In particular the extreme points of $\overline S _H^q(m,\alpha)$ are 
$\left\{{{h_u}}\right\}$ and $\left\{{{g_u}}\right\}$.\\
\end{theorem}
The following theorem gives the bounds for functions in 
$\overline S _H^q(m,\alpha,u)$ which yields a covering result for this class.\\
\begin{theorem}
Let $f \in \overline S _H^q(m,\alpha,u)$.Then for 
$\left| z \right| = r < 1$, we have 
\begin{equation*}
\left| {f(z)} \right| \le \left( {1  +\left| {{b_1}} \right|} \right)r + \frac{1}{{{2^n}}}\left( {1 - \left| {{b_1}} \right| - \alpha }\right){r^2},\quad \left| z \right| = r < 1
\end{equation*}
and
\begin{equation*}
\left| {f(z)} \right| \ge \left( {1 - \left| {{b_1}} \right|} \right)r - \frac{1}{{{2^n}}}\left( {1 - \left| {{b_1}} \right| - \alpha }\right){r^2},\quad  \left| z \right| = r < 1.
\end{equation*}
\end{theorem}
\textbf{Proof:}
Let $f \in \overline S _H^q(m,\alpha,u)$.Taking the absolute value of $f(z)$ ,we have\\
\begin{align*}
\left|f(z)\right|&\le(1+\left|b_1\right|)r+\sum\limits_{u \ge 2}^\infty  {\left( {\left| {{a_u}} \right| + \left| {{b_u}} \right|} \right)} {r^u}\\
&\le \left( {1 + \left| {{b_1}} \right|} \right)r + \sum\limits_{u \ge 2}^\infty  {\left( {\left| {{a_u}} \right| + \left| {{b_u}} \right|} \right)} {r^2}\\
&\le \left( {1 + \left| {{b_1}} \right|} \right)r + \frac{1}{{\left[ 2 \right]_q^m}}\sum\limits_{u \ge 2}^\infty  {{\left[ u \right]_q^m}\left( {\left| {{a_u}} \right| + \left| {{b_u}} \right|} \right)} {r^2}\\
&\le \left( {1 + \left| {{b_1}} \right|} \right)r + \frac{1}{{\left[ 2 \right]_q^m}}\left( {1 - \alpha  - \left| {{b_1}} \right|} \right){r^2}
\end{align*}
and 
\begin{align*}
\left|f(z)\right|&\ge(1-\left|b_1\right|)r-\sum\limits_{u \ge 2}^\infty  {\left( {\left| {{a_u}} \right| + \left| {{b_u}} \right|} \right)} {r^u}\\
&\ge \left( {1 - \left| {{b_1}} \right|} \right)r - \sum\limits_{u \ge 2}^\infty  {\left( {\left| {{a_u}} \right| + \left| {{b_u}} \right|} \right)} {r^2}\\
&\ge \left( {1 - \left| {{b_1}} \right|} \right)r - \frac{1}{{\left[ 2 \right]_q^m}}\sum\limits_{u \ge 2}^\infty  {{\left[ u \right]_q^m}\left( {\left| {{a_u}} \right| + \left| {{b_u}} \right|} \right)} {r^2}\\
&\ge \left( {1 - \left| {{b_1}} \right|} \right)r - \frac{1}{{\left[ 2 \right]_q^m}}\left( {1 - \alpha  - \left| {{b_1}} \right|} \right){r^2}
\end{align*}
The functions 
$z + \left| {{b_1}} \right|\overline z  + \frac{1}{{\left[ 2 \right]_q^m}}\left( {1 - \alpha  - \left| {{b_1}} \right|} \right){\overline z ^2}$ and $z - \left| {{b_1}} \right|z - \frac{1}{{\left[ 2 \right]_q^m}}\left( {1 - \alpha  - \left| {{b_1}} \right|} \right){z^2}$ for 
$\left| {{b_1}} \right| \le \left( {1 - \alpha } \right)$. 
\pagebreak


\begin{thebibliography}{17}
\bibitem{B1}
 J. Clunie as well as T. Sheil-Small, Harmonic univalent functionalities, Ann.Acad. Sci. Fenn.Ser. A I Mathematics. 9 (1984 ), 3-25.
\bibitem{B2}
J. M. Jahangiri, Accordant features starlike in the unit disk, J. Mahematics. Anal. Appl. 235 (1999), no. 2, 470-477.
\bibitem{B3}
J. M. Jahangiri, Harmonic univalent features determined next to q-calculus  drivers, Int. J.Mathematics. Anal. and also Appl. 5 (2018 ), no. 2, 39-43.
\bibitem{B4}
 K.K Dixit. and Porwal Saurabh, A subclass of harmonic univalent functions with positive coefficients, Tamkang J. Math., 41(3) (2010), 261-269. 
\bibitem{B5}
S.B. Joshi, and N.D.Sangle, New subclass of univalent functions defined by using generalised Salagean operator, J. Indones Mathematics Society (MIHMI)., 15 (2009), 79–89 .
\bibitem{B6}
S.D. Purohit  and  R.K. Raina, Certain subclasses of analytic functions
associated with fractional q-calculus operators, Math. Scand. 109 (2011), no. 1, 55-70.
\bibitem{B7} 
J.M.Jahangiri , Harmonic functions starlike in the unit disc, J. Math. Anal. Appl., 235 (1999), 470-477.
\bibitem{B8}
 N.D.Sangle,and  Y.P. Yadav, On a subclass of harmonic univalent functions defined By generalized derivative operator, IJMER, Vol.2 (2012), Issue 3, 562-569.
\bibitem{B9}
S.B. Joshi and N.D. Sangle, New subclass of Goodman-typep-valent harmonic functions, Filomat, 22(1) (2008),  193–204.
\bibitem{B10}
G.S. Salagean , Subclasses of univalent functions, Complex Analysis-Fifth Romanian Finish Seminar, Bucharest, 1(1983) 362-372.
\bibitem{B11}
H. Silverman , Harmonic univalent function with negative coefficients, J. Math. Anal. Appl., 220 (1998), 283-289.
\bibitem{B12}
H. Silverman  and Silvia E.M. , Subclasses of Harmonic univalent functions, New Zealand J. Math., 28(1999), 275-284.
\bibitem{B13}
Balvir Singh and Porwal Saurabh ,On A New Subclass of A Harmonic Univalent Functions, IJCRT, Vol.5(2017), Issue 4, Page 3465-3469.
\bibitem{B14} 
B. Ravindar, R.B.Sharma and N. Magesh, On Certain Subclass of Harmonic Univalent Functions Defined Q-Differential Operator, J.Mech.Cont.Math.Sci., Vol.14(2019), No.6, 45-53. 
\bibitem{15}
N. D. Sangle and G. M. Birajdar, Certain subclass of analytic function with negative coefficients defined by Catas operator, Indian Journal of Mathematics(IJM), 62(3) (2020),  335-353.
\end{thebibliography}
\end{document}